\documentclass[11pt]{amsart}
\numberwithin{equation}{section}
\newtheorem{theorem}{Theorem}
\newtheorem{proposition}{Proposition}

\newtheorem{definition}{Definition}
\numberwithin{theorem}{section} \numberwithin{lemma}{section}
\numberwithin{proposition}{section}
\numberwithin{definition}{section}
\def\R{\bf R}

\def\al{\aligned}
\def\eal{\endaligned}
\def\M{{\bf M}}
\def\be{\begin{equation}}
\def\ee{\end{equation}}
\def\lab{\label}

\begin{document}

\tracingpages 1
\title[ancient solutions]{\bf Heat kernel bounds, ancient $\kappa$ solutions and
the Poincar\'e conjecture }
\author{Qi S. Zhang}
\address{ Department of
Mathematics, University of California, Riverside, CA 92521, USA }

\begin{abstract}
We establish certain Gaussian type upper bound for the heat kernel
of the conjugate heat equation associated with
3 dimensional ancient $\kappa$ solutions to the Ricci flow.

As an application, using the $W$ entropy associated with the heat kernel,
we give a different and much shorter proof of Perelman's
classification of backward limits of these ancient solutions. The method is partly motivated by \cite{Cx:1} and \cite{S:1}.  The
current paper or \cite{CL:1} combined with \cite{ChZ:1} and \cite{Z:2}  lead to a
simplified proof of the Poincar\'e conjecture without using reduced
distance and reduced volume.

\end{abstract}
\maketitle \tableofcontents
\section{Introduction}

The main goal of the paper is to establish certain Gaussian type upper bound
for the heat kernel (fundamental solutions) of the conjugate heat equation
 associated with
3 dimensional ancient $\kappa$ solutions to the Ricci flow. Heat
kernel estimates have been an active area of research.  
 When coupled
with Ricci flow, various estimates can be found in \cite{G:1},
\cite{P:1} Section 9, \cite{Ni:1}, and \cite{Z:1}.  For example, in Section 9 of \cite{P:1}
Perelman proved a lower bound for the fundamental solution of the conjugate heat 
equation for general Ricci flow. So far an upper bound corresponding to this lower bound has been missing.  Our result is a progress in this direction when the Ricci 
flow is a 3 dimensional ancient $\kappa$ solution.

One motivation of the work is that it induces a simpler proof of  the Poincar\'e conjecture.  The most difficult analytical parts of the proof can now be treated by one unifying theme:  Perelman's $W$ entropy and related 
(log) Sobolev inequalities and heat kernel estimates.  Let us explain the point in 
more detail.
From Perelman's original papers \cite{P:1}, \cite{P:2}, \cite{P:3}
and the works by Cao and Zhu \cite{CZ:1}, Kleiner and Lott
\cite{KL:1} and Morgan and Tian \cite{MT:1}, and Tao \cite{T:2},
\cite{T:1},
 it is clear that
the bulk of the proof of the Poincar\'e conjecture is consisted of
two items. One is the proof of local non-collapsing with or without
surgeries, and the other is the classification of backward limits of
ancient $\kappa$ solutions. After these are done, one can show that
regions where the Ricci flow is close to forming singularity have
simple topological structure, i.e. canonical neighborhoods. Then one
proceeds to prove that the singular region can be removed by finite
number of surgeries in finite time. When the initial manifold is
simply connected, the Ricci flow becomes extinct in finite time [P3]
(see also \cite{CM:1}).  Thus the manifold is diffeomorphic to ${\bf
S}^3$, as conjectured by Poincar\'e.

Besides the results and techniques by R. Hamilton, the main new
tools Perelman used in carrying out the proof are several monotone
quantities along Ricci flow. These include the $W$ entropy, reduced
volume and the associated reduced distance. In \cite{P:1}, Perelman
first used his $W$ entropy to prove local non-collapsing for smooth
Ricci flows.  However he then turned to the reduced
volume (distance) to prove the classification and non-collapsing
with surgeries. The $W$ entropy is not used anymore. The reduced
distance,  not being smooth or positive in general,
 is one of the causes of the complexity of the original proof.

It turns out that the $W$ entropy is just the formula in a log Sobolev inequality (c.f. 
\cite{Gr:1} in the fixed metric case) and
the monotonicity of the $W$ entropy implies certain uniform Sobolev inequalities
along the Ricci flow. Using this idea and being inspired by the last section of [P2]
 and [KL], we proved in [Z2]  a stronger local non-collapsing result for Ricci flow
 with surgeries. The proof, without using reduced distance or volume,
 is short and seems more accessible.
 It also strengthens and clarifies the original result by doing analysis at one
time level each time, thus avoiding the complication associated with
surgeries. In the wake of this development, it would be desirable
that the classification
 mentioned above can also be done by using the $W$ entropy alone.
 Such a view was also expressed in \cite{T:1} e.g.

As one application of the main result of the paper, using the $W$ entropy associated with the heat
kernel,
 we give a different and much shorter proof of Perelman's
classification of backward limits of these ancient solutions. Thus,
the current paper together with \cite{Z:2} and \cite{ChZ:1}(see
explanation 4 paragraphs below) lead to a simplified proof of the
Poincar\'e conjecture. Of course we still follow the framework by
Perelman. However, much of the highly intensive analysis involving
reduced distance and volume is now replaced by the study of the $W$
entropy and the related uniform Sobolev inequalities and heat kernel
estimates. Sobolev inequalities and heat kernels are familiar to
many mathematicians. Therefore the current proof is more accessible
to a wider audience. Besides, due to the relative simplicity, we
hope the current technique can lead to better understanding of other
problems for Ricci flow.

We should mention that the reduced distance and volume are still needed for the
proof of the geometrization conjecture. Specifically, they are needed,
but only in the proof of Perelman's no local collapsing Theorem II with surgeries.

Let us outline the proof. In the next section we prove Theorem 1.1
concerning the bounds for the heat kernel of the conjugate heat
equation.  The proof follows the framework in section 5 of [Z1].
There an upper bound in the case of Ricci flow with nonnegative
Ricci curvature was given. In the current situation, the  ancient
$\kappa$ solutions provide better control on curvature and volume.
These allow us to find a better Gaussian upper bound for the heat
kernel.  These bounds can be regarded as generalization of the heat
kernel bounds of Li and Yau \cite{LY:1} in the fixed metric case.

Using this heat kernel bound, in Section 3 we show that the $W$
entropy associated with the heat kernel is uniformly bounded from
below after certain scaling. After this done, we use Perelman's
monotonicity formula for the $W$ entropy to prove the backward limit
is a shrinking gradient Ricci soliton. This part of the arguments
resembles that in the paper \cite{Cx:1} and \cite{S:1}  where
forward convergence results for normalized Ricci flow were proven.

Finally one needs to prove universal non-collapsing for ancient
$\kappa$ solutions without reduced distance or volume.  But this is
already done in \cite{ChZ:1} , even in certain more general 4
dimensional situation. We will just describe their proof.

Now let us introduce the definitions and notations in order to
present our result precisely.  ${\bf M}$ denotes a complete compact,
or noncompact Riemannian manifold, unless stated otherwise; $g,
R_{ij}$ (or $Ric$) will be the metric and Ricci curvature; $\nabla$,
$\Delta$ the corresponding gradient and Laplace-Beltrami operator;
$c$ with or without index denote generic positive constant that may
change from line to line. If the metric $g(t)$ evolves with time,
then $d(x, y, t)$ will denote the corresponding distance function;
$dg(x, t)$  or $dg(t)$ denote the volume element under $g(t)$; We
will use $B(x, r; t)$ to denote the geodesic ball centered at $x$
with radius $r$ under the metric $g(t)$; $|B(x, r; t)|_s$ to denote
the volume of $B(x, r; t)$ under the metric $g(s)$.
 We will still use $\nabla$, $\Delta$  to denote the
corresponding gradient and Laplace-Beltrami operator for $g(t)$,
without mentioning the time $t$, when no confusion arises.

We use the following concept of ancient $\kappa$ solutions according to
Perelman.

\begin{definition}
\lab{defgujie}
A solution to the Ricci flow $\partial_t g = - 2 Ric$ is an ancient $\kappa$ solution
if it satisfies the following properties.

1. It is complete (compact or noncompact) and defined on an ancient time interval
$(-\infty, T_0]$,  $T_0 \ge 0$.

2. It has nonnegative curvature operator and bounded curvature at each time level.

3.  It is $\kappa$ noncollapsed on all scales for some positive constant $\kappa$.
i.e.

Suppose that $x_0 \in {\bf M}$, $t_0 \in (-\infty, T_0]$.
 Let $P(x_0, t_0, r, -r^2)$ be the  parabolic ball
 \[
 \{(x, t) \ | \
 d(x, x_0, t) < r, \quad t_0-r^2 < t < t_0 \ \}.
 \]
  Then ${\bf M}$ is $\kappa$ non-collapsed at $(x_0, t_0)$ at
scale $r$ if $|Rm| \le r^{-2}$ on $P(x_0, t_0, r, -r^2)$ and
$vol(B(x_0, t_0, r)) \ge  \kappa r^3$.
\end{definition}
\medskip

For convenience, we take the final time $T_0$ of the ancient
solution to be $0$ throughout the paper.  The conjugate heat
equation is
\be
\lab{eqconj}
 \Delta u - R u - \partial_\tau u=0.
\ee     Here and always $\tau = - t$.  $\Delta$ and $R$ are the
Laplace-Beltrami operator and the scalar curvature with respect to
$g(t)$. This equation, coupled with the initial value $u_{\tau=0} =
u_0$ is well posed if $\M$ is compact or the curvature is bounded,
and if $u_0$ is bounded \cite{G:1}..

We use $G=G(x, \tau; x_0, \tau_0)$ to denote the heat kernel
(fundamental solution) of (\ref{eqconj}).  Here $\tau>\tau_0$ and
$x, x_0 \in \M$. Existence of $G$ was established in \cite{G:1}.
The main result of the
paper is

\begin{theorem}
\label{th1.1} (i). Let $(\M, g(t))$ be a $n$ dimensional ancient
$\kappa$ solution of the Ricci flow.  Suppose also that $R(x, t) \le
\frac{D_0}{1+|t|}$ for some $D_0>0$ and for $t \in [-T, 0]$. Here
$T$ is any positive number or $T=\infty$.  Then exist positive
numbers $a$ and $b$ depending only on $n$,  $\kappa$ and $D_0$ such
that the following holds.

For all $x, x_0 \in \M$,
\[
G(x, \tau; x_0, \tau_0) \le \frac{a}{|B(x, \sqrt{\tau-\tau_0}, t_0)|_{t_0}}
e^{-b d^2(x, x_0, t_0)/(\tau-\tau_0)},
\]where $\tau = -t$, $\tau_0=-t_0$, $\tau>\tau_0$ and $ t \in [-T, 0]$.

(ii). In particular, if $R(x, t) \le \frac{D_0}{1+|t|}$ for all $ t \le 0$, namely $(\M, g(t))$ is a
Type I ancient solution, there exist positive numbers $a_1$
and $b_1$ depending only on $\kappa$ and $D_0$ such that the following holds.
For all $x, x_0 \in \M$,  and all $\tau = - t >0$,
\[
\frac{1}{a_1 \tau^{n/2}}
e^{- d^2(x, x_0, t)/(b_1 \tau)} \le G(x, \tau; x_0, \tau/2) \le \frac{a_1}{\tau^{n/2}}
e^{-b_1 d^2(x, x_0, t)/\tau}.
\]

\end{theorem}

{\it Remark.} The full Gaussian lower bound in part (ii) of the
theorem is not needed for the application in Section 3. One only
needs the lower bound for one point in the ball $B(x_0, \sqrt{b|t|},
t)$ for some $b>1$, which is a simple consequence of the upper bound.

The Gaussian upper and lower bounds seem to be of interest that is
independent of the Poicar\'e conjecture.  For instance, Perelman \cite{P:1} used 
heat kernel bounds to prove his pseudo locality theorem.  In Section 9 of the same paper, a lower bound for the heat kernel was proven.  However the upper bound 
is missing.
 In this sense, this paper
is not just a reproof of a known result.

\section{Proof of Theorem \ref{th1.1}:  the heat kernel  bounds}

We divide the proof into three steps. The first two are for part (i)
of the theorem.  We always assume that all the time variables
involved are not smaller than $-T$, so that the condition $R(\cdot.
t) \le \frac{D_0}{1+|t|}$ holds.  As mentioned in the introduction,
the proof follows the framework of Theorem 5.2 in \cite{Z:1}, where
certain upper bound for $G$ under Ricci flow with nonnegative Ricci
curvature was derived. Comparing with that case, we have two new
ingredients coming from ancient $\kappa$ solutions. One is the
non-collapsing condition on all scales. The other is the bound on
the scalar curvature.  These allow us to prove a better bound.
During the proof, there will be overlaps with \cite{Z:1}. They are
here so that the paper is self contained. Without loss of generality
we assume $\tau_0=0$ in $G(x, \tau; x_0, \tau_0)$.   It is
convenient to work with the reversed time $\tau$. Note that the
Ricci flow is a backward flow with respect to $\tau$ and the
conjugate heat equation is a forward heat equation with a potential
term.

{\bf Step 1.}

Since $Ricci \ge 0$, it is well known (see Theorem 3.7 \cite{Heb:1}
e.g.) the following Sobolev inequality holds: Let $B(x, r, t)$ be a
proper subdomain for $(\M, g(t))$.  For all $v \in W^{1, 2}(B(x, r,
t))$, there exists $c_n>0$ depending only on the dimension $n$ such
that \be \lab{SobRic>0} \bigg{(}  \int v^{2n/(n-2)} dg(t)
\bigg{)}^{(n-2)/n} \le \frac{c_n  r^2}{|B(x,  r, t)|^{2/n}_t}  \int
\left[ | \nabla v |^2  + r^{-2} v^2 \right] dg(t). \ee

For our purpose, we only need to take $r = c \sqrt{|t|}$, for $c<1$.
By the assumption that $R(x, t) \le \frac{D_0}{1+|t|}$ and the
$\kappa$ non-collapsing property, we have
\[
| B(x, \sqrt{|t|},  t)|_t \ge \kappa D^{-n}_0 |t|^{n/2}.
\]Therefore the above Sobolev inequality becomes
\be
\lab{SobT}
\bigg{(}  \int v^{2n/(n-2)} dg(t) \bigg{)}^{(n-2)/n} \le
\frac{c_n  D^2_0}{\kappa^{2/n}}  \int \left[ | \nabla
v |^2  + |t|^{-1} v^2 \right] dg(t)
\ee for all $v \in W^{1, 2}(B(x, \sqrt{|t|},  t))$.

Before moving forward, we would like to clarify a technical point in the definition
of Perelman's $\kappa$ non-collapsing as given in Definition \ref{defgujie}.  
The issue is whether the metric balls $B(x, r, t)$ in the definition are required to be a proper subdomain
of the manifold $\M$.  When $\M$ is noncompact, $B(x, r, t)$ is always a proper 
subdomain so this issue is mute.  Now one assumes that $\M$ is compact.  Without requiring $B(x, r, t)$ being 
a proper subdomain, if $r$ is larger than the diameter of $\M$, then $B(x, r, t)$ is the whole manifold. In this case $|B(x, r, t)|_t$ can not be greater than $\kappa r^n$ for large $r$. So to be $\kappa$ non-collapsed,  at some point in the parabolic ball $|Rm|$ is 
greater than $1/r^2$.  In other words, if $|Rm| \le 1/r^2$ in the parabolic ball, then
the volume of the manifold is at least $\kappa r^n$.  If the Ricci curvature is nonnegative, then by standard volume comparison theorem,  the diameter of the 
manifold at time $t$ is at least $c r$. 

In this paper, we take this explanation for Perelman's $\kappa$ non-collapsing, i.e $B(x, r, t)$ in the definition of $\kappa$ ancient solutions is not required to be a proper 
subdomain. This seems to be the prevailing view in the literature.  That is why 
  the Sobolev imbedding \ref{SobT} holds without requiring that $B(x, \sqrt{|t|}, t)$
  is a proper subdomain of $\M$.   A natural question is: what happens when 
  $B(x, r, t)$ is implicitly assumed as a proper subdomain in the definition of 
  $\kappa$ solutions?  Then we have to make this extra assumption throughout. However
 either way does not
affect the application for the Poincar\'e conjecture in the next section. The reason is 
compact ancient solutions are already taken care of. See the beginning of the proof of Theorem \ref{thbacklimit}.

Next we show that, under the assumptions of the theorem, $(\M, g(t))$
possess a space time doubling property:  the distance between two points at times
$t_1$ and $t_2$ are comparable if $t_1$ and $t_2$ are comparable. The proof is very
simple.  Given $x_1, x_2 \in \M$, let $\bf r$ be a shortest geodesic connecting the two.
Then
\[
\partial_t d(x_1, x_2, t) = - \int_{\bf r} Ric (\partial_r, \partial_r) ds.
\]Since the sectional curvature is nonnegative, it holds
\[
| Ric(x, t) | \le c R(x, t) \le \frac{c D_0}{1+|t|}.
\]Therefore
\[
 -\frac{c D_0}{1+|t|} d(x_1, x_2, t) \le \partial_t d(x_1, x_2, t) \le 0.
\]After integration, we arrive at:
\be
\lab{distdouble}
\left( |t_1/|t_2| \right)^{c D_0} \le d(x_1, x_2, t_1)/d(x_1, x_2, t_2) \le 1
\ee  for all $t_2<t_1<0$.
Note that the above inequality is of local nature. If the distance is not smooth, then
one can just shift one point, say $x_1$,  slightly and then obtain the same integral
inequality by taking limits.

Similarly, we have
\[
0 \ge \partial_t \int_{B(x, \sqrt{|t_1|}, t_1)} dg(t) =
- \int_{B(x, \sqrt{|t_1|}, t_1)} R(y, t) dg(t) \ge - \frac{ D_0}{1+|t|}
 \int_{B(x, \sqrt{|t_1|}, t_1)} dg(t).
\]Upon integration, we know that the volume of the balls
\be
\lab{voldouble}
|B(x,  \sqrt{|t_3|}, t_4)|_{t_5}
\ee  are all comparable for $t_3, t_4, t_5 \in [t_2, t_1]$, provided that
$t_1$ and $t_2$ are comparable.

Let $u$ be a positive solution to (\ref{eqconj}) in the region
\[
Q_{\sigma r}(x, \tau) \equiv \{ (y, s) \ | \ y \in {\bf M}, \tau-(\sigma
r)^2 \le s \le \tau, \ d(y, x, -s) \le \sigma r \}.
\]Here $r=\sqrt{|t|}/8>0, 2 \ge \sigma \ge 1$.
Given any $p \ge 1$, it is clear that
\be
\lab{equp}
\Delta u^p - p R u^p- \partial_\tau u^p \ge 0.
\ee

Let $\phi: [0, \infty) \to [0, 1]$ be a smooth function such that
$|\phi'| \le 2/((\sigma-1) r)$, $\phi' \le 0$, $\phi(\rho) =1$ when
$0 \le \rho \le r$, $\phi(\rho)=0$ when $\rho \ge \sigma r$. Let
$\eta: [0, \infty) \to [0, 1]$ be a smooth function such that
$|\eta'| \le 2/((\sigma-1) r)^2$, $\eta' \ge 0$, $\eta \ge 0$,
$\phi(s) =1$ when $\tau-r^2 \le s \le \tau$, $\phi(s)=0$ when $s \le
\tau- (\sigma r)^2$.  Define a cut-off function $\psi = \phi(d(x, y,
-s)) \eta(s)$.

Writing $w= u^p$ and using $w \psi^2$ as a test function on (\ref{equp}),
we deduce
\be
\lab{dwdw}
\int \nabla (w \psi^2) \nabla w dg(y, -s)ds + p \int R w^2 \psi^2
dg(y, -s)ds \le - \int (\partial_s w) w \psi^2  dg(y, -s)ds.
\ee  By direct calculation
\[
\int \nabla (w \psi^2) \nabla w dg(y, -s)ds = \int  | \nabla (w
\psi)|^2 dg(y, -s)ds - \int  |\nabla \psi|^2 w^2 dg(y, -s)ds.
\]Next we estimate the righthand side of (\ref{dwdw}).
\[
\al
- \int (\partial_s w)& w \psi^2  dg(y, -s)ds\\
& = \int w^2 \psi
\partial_s \psi dg(y, -s)ds + \frac 1 2 \int (w \psi)^2 R dg(y, -s)ds
- \frac{1}{2} \int (w \psi)^2 dg(y, -\tau).
\eal
\]Observe that
\[
 \partial_s \psi = \eta(s)   \phi'(d(y, x, -s)) \partial_s d(y, x, -s)
+ \phi(d(y, x, -s)) \eta'(s) \le \phi(d(y, x, -s)) \eta'(s).
\]This is so because $\phi' \le 0$ and $\partial_s d(y, x, -s) \ge 0$
under the Ricci flow with nonnegative Ricci curvature.
Hence
\be
\lab{dwdw2}
\al
 - \int &(\partial_s w) w \psi^2  dg(y, -s)ds \\
 &\le \int w^2 \psi \phi(d(y, x, -s))  \eta'(s)  dg(y, -s)ds +
 \frac 1 2 \int (w \psi)^2 R dg(y, -s)ds
- \frac{1}{2} \int (w \psi)^2 dg(y, -\tau). \eal
\ee  Combing (\ref{dwdw}) with (\ref{dwdw2}),
we obtain, in view of $p \ge 1$ and $R \ge
0$,
\be
\lab{dwdw3}
\int  | \nabla (w \psi)|^2 dg(y, -s)ds + \frac{1}{2} \int (w
\psi)^2 dg(y, -\tau) \le \frac{c}{(\sigma-1)^2 r^2} \int_{Q_{\sigma r
(x, \tau)}} w^2 dg(y, -s)ds.
\ee   By H\"older's inequality
\be
\lab{2/n}
\int (\psi w)^{2(1+(2/n)} dg(y, -s) \le \bigg{(}   \int (\psi
w)^{2n/(n-2))} dg(y, -s) \bigg{)}^{(n-2)/n} \bigg{(} \int (\psi
w)^2 dg(y, -s) \bigg{)}^{2/n}.
\ee
\medskip

  By the $\kappa$ non-collapsing assumption,
$|B(x, \sqrt{|t|}, t)|_{t} \ge \kappa c_2 r^n$.  Since $\M$ has
nonnegative Ricci curvature, the diameter of ${\M}$ at time $t$ is a
least a constant multiple of $c \sqrt{|t|}$ for some $c=c_n>0$.
 Recall that $r =
\sqrt{|t|}/8$. Therefore by the distance doubling property
(\ref{distdouble}),   $B(x, \sigma r, -s)$ is a proper sub-domain of
${\bf M}$, $s \in [\tau-(\sigma r)^2, \tau]$.  Here we just take the
number $8$ for simplicity. If it is not large enough, we just
replace it by a sufficiently large number $D$ and consider
$r=\sqrt{|t|}/D$ instead. By the Sobolev inequality (\ref{SobT}), it
holds
\[
\bigg{(}  \int (\psi w)^{2n/(n-2)} dg(y, -s) \bigg{)}^{(n-2)/n} \le
c(\kappa, D_0)
\int [| \nabla (\psi w) |^2  + r^{-2} (\psi w)^2] dg(y, -s),
\]for $s \in [t-(\sigma r)^2, t]$. Substituting  this and (\ref{dwdw3}) to (\ref{2/n}), we arrive at the estimate
\[
\int_{Q_r(x, \tau)} w^{2 \theta} dg(y, -s) ds \le c(\kappa, D_0) \bigg{(} \frac{1}{(\sigma-1)^2
r^2} \int_{Q_{\sigma r}(x, \tau)} w^2 dg(y, -s) ds \bigg{)}^{\theta},
\]with $\theta = 1+(2/n)$. Now we apply the above inequality repeatedly with
the parameters $\sigma_0=2, \sigma_i=2- \Sigma^i_{j=1} 2^{-j}$ and
$p=\theta^i$. This shows a $L^2$ mean value inequality \be
\lab{L2mvi} \sup_{Q_{r/2}(x, \tau)} u^2 \le \frac{c(\kappa,
D_0)}{r^{n+2} } \int_{Q_r(x, \tau)} u^2 dg(y, -s)ds. \ee

This inequality clearly also holds if one replaces $r$ by any
positive number $r' <r$ since $|B(x, r', t)| \ge k c_n |B(x, r, t)|
(r'/r)^n \ge c r'^n$ by the doubling condition for manifolds with
nonnegative Ricci curvature. Then one can just rerun the above
Moser's iteration.

  From
here, by a generic trick of Li and Schoen \cite{LS:1}, applicable
here since it
 uses
only the doubling property of the metric balls, we arrive at the
$L^1$ mean value inequality
\[
\sup_{Q_{r/2}(x, \tau)} u \le  \frac{c(\kappa, D_0)}{r^{n+2} }
\int_{Q_r(x, \tau)} u dg(z, -s)ds.
\]We remark that the doubling constant is uniform since the metrics have nonnegative
Ricci curvature.

Now we take  $u(x, \tau) = G(x, \tau; x_0, 0)$. Note that $\int_{\bf
M} u(z, s) dg(z, -s) =1$ and $r=\sqrt{|t|}$. \be \lab{Gdbound} G(x,
\tau; x_0, 0) \le \frac{c(\kappa, D_0)}{|t|^{n/2}}. \ee

\medskip

{\bf step 2.}  proof of the Gaussian upper bound.

 We begin by using a
modified version of the exponential weight method due to Davies
\cite{Da:1}. Pick a point $x_0 \in {\bf M}$, a number $\lambda<0$ and a
function $f \in C^{\infty}_0({\bf M}, g(0))$. Consider the functions $F$ and
$u$ defined by
\be
\lab{f&u}
F(x, \tau) \equiv  e^{\lambda d(x, x_0, t)} u(x, \tau) \equiv e^{\lambda
d(x, x_0, t)} \int G(x, \tau; y, 0) e^{-\lambda d(y, x_0, 0)} f(y)
dg(y, 0).
\ee   Here and always $\tau = -t$.  It is clear that $u$ is a solution of (\ref{eqconj}).  By direct
computation,
\[
\al
\partial_\tau &\int F^2(x, \tau) dg(x, t) = \partial_\tau \int
e^{2\lambda d(x, x_0, t)} u^2(x, \tau)dg(x, t)\\
&=  2\lambda \int e^{2\lambda d(x, x_0, t)} \partial_\tau d(x, x_0, t)
u^2(x, \tau)dg(x, t) + \int e^{2\lambda d(x, x_0, t)} u^2(x, \tau) R(x,
t) dg(x, t)\\
&\qquad + 2 \int e^{2\lambda d(x, x_0, t)}[\Delta u - R(x, t) u(x,
\tau)] u(x, \tau)dg(x, t).
 \eal
\]By the assumption that $Ricci \ge 0$ and $\lambda<0$, the above
shows
\[
\partial_\tau  \int F^2(x, \tau) dg(x, t)
\le 2  \int e^{2\lambda d(x, x_0, t)} u \Delta u (x, \tau) dg(x, t).
\]Using integration by parts, we turn the above inequality into
\[
\al
\partial_\tau &\int F^2(x, \tau) dg(x, t)\\
&\le -4  \lambda \int e^{2\lambda d(x, x_0, t)} u \nabla d(x, x_0,
t) \nabla u
 dg(x, t) -2 \int e^{2\lambda d(x, x_0, t)} |\nabla u|^2 dg(x,
t). \eal
\]Observe also
\[
\al \int &|\nabla F(x, \tau)|^2 dg(x, t) = \int |\nabla (e^{\lambda
d(x, x_0, t)} u(x, \tau))|^2 dg(x, t)\\
&=\int e^{2\lambda d(x, x_0, t)} |\nabla u|^2 dg(x, t) + 2 \lambda
\int e^{2\lambda d(x, x_0, t)} u \nabla d(x, x_0, t) \nabla u
 dg(x, t)\\
 &\qquad + \lambda^2 \int e^{2\lambda d(x, x_0, t)} |\nabla d|^2 u^2
 dg(x, t).
 \eal
 \]Combining the last two expressions, we deduce
 \[
\partial_\tau  \int F^2(x, \tau) dg(x, t)
\le - 2 \int |\nabla F(x, \tau)|^2 dg(x, t) + \lambda^2 \int
e^{2\lambda d(x, x_0, t)} |\nabla d|^2 u^2
 dg(x, t).
 \]By the definition of $F$ and $u$, this shows
 \[
\partial_\tau  \int F^2(x, \tau) dg(x, t) \le \lambda^2 \int F(x, \tau)^2
 dg(x, t).
\]Upon integration, we derive the following $L^2$ estimate
\be
\lab{F2t}
\int F^2(x, \tau) dg(x, t) \le e^{\lambda^2 \tau} \int F^2(x, 0) dg(x,
0) = e^{\lambda^2 \tau} \int f(x)^2 dg(x, 0).
\ee

Recall that $u$ is a solution to (\ref{eqconj}). Therefore, by the mean value
inequality (\ref{L2mvi}), the following holds
\[
u(x, \tau)^2 \le \frac{c(\kappa, D_0)}{\tau^{1+n/2}} \int^{\tau}_{\tau/2}
\int_{B(x, \sqrt{|t|/2}, -s)} u^2(z, s) dg(z, -s) ds.
\]By the definition  of $F$ and $u$, it follows that
\[
u(x, \tau)^2 \le  \frac{c(\kappa, D_0)}{\tau^{1+n/2}} \int^{\tau}_{\tau/2}
\int_{B(x, \sqrt{|t|/2}, -s)}  e^{-2 \lambda d(z, x_0, -s)} F^2(z,
s) dg(z, -s) ds.
\]In particular, this holds for $x=x_0$. In this case, for $z \in
B(x_0, \sqrt{|t|/2}, -s)$, there holds $d(z, x_0, -s) \le
\sqrt{|t|/2}.$ Therefore, by the assumption that $\lambda<0$,
\[
u(x_0, \tau)^2 \le  \frac{c(\kappa, D_0)}{\tau^{1+n/2}} e^{- \lambda
\sqrt{2|t|}} \int^{\tau}_{\tau/2} \int_{B(x_0, \sqrt{|t|/2}, -s)}
F^2(z, s) dg(z, -s) ds.
\]This combined with (\ref{F2t}) shows that
\[
u(x_0, \tau)^2 \le  \frac{c(\kappa, D_0)}{\tau^{n/2}} e^{\lambda^2 \tau
- \lambda \sqrt{2|t|}}  \int f(y)^2 dg(y, 0).
\]i.e.
\[
\bigg{(} \int G(x_0, \tau; z, 0) e^{-\lambda d(z, x_0, 0)} f(z) dg(z,
0) \bigg{)}^2 \le
 \frac{c(\kappa, D_0)}{\tau^{n/2}} e^{\lambda^2 \tau
- \lambda \sqrt{2|t|}}
 \int f(y)^2 dg(y, 0).
\]Now, we fix $y_0$ such that $d(y_0, x_0, 0)^2 \ge 4  t$. Then it is clear that, by $\lambda<0$ and the triangle
inequality,
\[
-\lambda d(z, x_0, 0) \ge - \frac{\lambda}{2} d(x_0, y_0, 0)
\]when $d(z, y_0, 0) \le \sqrt{|t|}$.
In this case, the above integral inequality implies
\[
\bigg{(} \int_{B(y_0, \sqrt{|t|}, 0)}  G(x_0, \tau; z, 0) f(z) dg(z, 0)
\bigg{)}^2 \le \frac{c(\kappa, D_0) e^{ \lambda d(x_0, y_0, 0) + \lambda^2
\tau- \lambda \sqrt{2|t|}}}{\tau^{n/2}}  \int f(y)^2
dg(y, 0).
\]Note that this inequality hold for all $-T \le t<0$ and $\lambda<0$. For an arbitrarily  fixed $t \in [-T, 0]$, 
 we take
\[
\lambda = - \frac{d(x_0, y_0, 0)}{\beta \tau}
\]with $\beta>0$  sufficiently large.
Since $f$ is arbitrary, this shows, for some $b>0$,
\[
\int_{B(y_0, \sqrt{|t|}, 0)}  G^2(x_0, \tau; z, 0) dg(z, 0) \le
\frac{c(\kappa, D_0) e^{-b d(x_0, y_0, 0)^2/\tau}}{\tau^{n/2}}.
\]Hence, there exists $z_0 \in B(y_0, \sqrt{|t|}, 0)$ such that
\[
G^2(x_0, \tau; z_0, 0) \le \frac{c(\kappa, D_0) }{
\tau^{n/2}  \ |B(x_0, \sqrt{|t|}, 0)|_0}   e^{-b d(x_0, y_0, 0)^2/\tau}.
\]

In order to get the upper bound for all points, let us consider the function
\[
v=v(z, l) \equiv G(x_0, \tau; z, l).
\]This is a
solution to the conjugate of the conjugate equation (\ref{eqconj}). i.e.
 \[ \Delta_z G(x, \tau; z; l) + \partial_l G(x, \tau; z, l) =0, \quad
 \partial_l g = 2 Ric.
 \]

 Therefore, we can use
Theorem 3.3 in \cite{Z:1},  after a reversal in time.  Note this
theorem was stated only for compact manifolds. However, as remarked
there, it is valid in he noncompact case whenever the maximum
principle for the heat equation holds. Since the proof is quite
short, we will present it in the appendix.  It is just a simple
generalization of Hamilton's first result in \cite{H:1} to the Ricci
flow case. Consequently, for $\delta>0, C>0$,
\[
G(x_0, \tau; y_0, 0) \le C G^{1/(1+\delta)} (x_0, \tau, z_0, 0)
M^{\delta/(1+\delta)},
\]where $M =\sup_{\M \times [0, \tau/2]} G(x_0, \tau, \cdot, \cdot)$.
By Step1, there exists a constant $c(\kappa, D_0)>0$, such that
\[
M \le \frac{c(\kappa, D_0)}{\tau^{n/2}}.
\]

Consequently
\[
G^2(x_0, \tau; y_0, 0) \le \frac{c(\kappa, D_0)}{\tau^{n/2} |B(x_0, \sqrt{|t|}, 0)|_0}
  \   e^{-b d(x_0, y_0, 0)^2/t} \le
   \frac{c(\kappa, D_0)}{ |B(x_0, \sqrt{|t|}, 0)|^2_0}
  \   e^{-b \  d(x_0, y_0, 0)^2/t}.
\]The last step holds since the Ricci curvature is nonnegative.

Since $x_0$ and $y_0$ are arbitrary, the proof of part (i) is done.

\medskip

{\bf step 3}

In this step, we prove the upper and lower bound for $G(x, \tau;
x_0, \tau/2)$ in the case of type I ancient solution. The upper
bound is already proven in view the distance and  volume  comparison
result (\ref{distdouble}), (\ref{voldouble}) and the fact that
$|B(x, \sqrt{|t|}, t)|_t \ge c(\kappa, D_0) |t|^{n/2}$. So we just
need to prove the lower bound.

For a number $\beta>0$ to be fixed later, the upper bound implies
\[
\al
\int_{B(x_0, \sqrt{ \beta |t|}, t)}& G^2(x, \tau; x_0, \tau/2) dg(x, t)\\
&\ge \frac{1}{|B(x_0, \sqrt{ \beta |t|}, t)|_t}
\left(  \int_{B(x_0, \sqrt{ \beta |t|}, t)} G(x, \tau; x_0, \tau/2) dg(x, t) \right) ^2\\
&=\frac{1}{|B(x_0, \sqrt{ \beta |t|}, t)|_t}
\left( 1- \int_{B(x_0, \sqrt{ \beta |t|}, t)^c} G(x, \tau; x_0, \tau/2) dg(x, t) \right) ^2\\
&\ge \frac{1}{|B(x_0, \sqrt{ \beta |t|}, t)|_t}
\left( 1- \int_{B(x_0, \sqrt{ \beta |t|}, t)^c} \frac{c(\kappa, D_0)}{ \tau^{n/2}}
  \   e^{-b \ d(x_0, y_0, t)^2/t} dg(x, t) \right) ^2
\eal
\]Since the Ricci curvature is nonnegative, one can use the volume doubling property
to compute that
\[
\int_{B(x_0, \sqrt{ \beta |t|}, t)^c} \frac{c(\kappa, D_0)}{ \tau^{n/2}}
  \   e^{-b \ d(x_0, y_0, t)^2/t} dg(x, t) \le 1/2
\]provided that $\beta$ is sufficiently large.
Here we stress that all constants are independent of $t$. Since
$|B(x_0, \sqrt{ \beta |t|}, t)|_t \le c_n(\beta |t|)^{n/2}$ by
standard volume comparison theorem, this shows
\[
\int_{B(x_0, \sqrt{ \beta |t|}, t)} G^2(x, \tau; x_0, \tau/2) dg(x, t) \ge
\frac{c(\kappa, D_0)}{ |t|^{n/2}}.
\]Hence there exists $x_1 \in B(x_0, \sqrt{ \beta |t|}, t)$ such that
\[
G(x_1, \tau; x_0, \tau/2) \ge \frac{c(\kappa, D_0)}{ |t|^{n/2}}.
\]For applications in Section 3, this lower bound is already
sufficient.

An inspection of the proof shows that actually for any $\lambda \in
[3/4, 4]$, it holds, for some $x_\lambda \in B(x_0, \sqrt{\beta
|t|}, t)$,
\[
G(x_\lambda, \lambda \tau; x_0, \tau/2) \ge \frac{c(\kappa, D_0)}{
|t|^{n/2}}.
\]It is well known that such a lower bound implies the full Gaussian lower bound if 
one has a suitable Harnack inequality. Such Harnack inequality already exists.  For the heat kernel, it is in Section 9 of
\cite{P:1}. For all positive solutions it is  in Corollary 2.1 (a) in \cite{KZ:1} and
\cite{CH:1}).  Applying Corollary 2.1 (a) in \cite{KZ:1}, we get
\[
G(x_{3/4}, \frac{3}{4} \tau; x_0, \tau/2)
 \leq
G(x,  \tau; x_0, \tau/2)
 \left(
            \frac{\tau }{\tau 3/4}
       \right)^n
 \exp
  {
   \frac
        {\int_{0}^{1}      [\,4|\gamma^{\prime}(s)|^2+(\tau/4)^2\,R\,]\,    ds}
        {2(\tau/4)}},
 \]where $\gamma$ is a smooth curve on $\M$ such that $\gamma(0)=x_{3/4}$
 and $\gamma(1)=x$. Also $|\gamma^{\prime}(s)|^2
 =g_{-l}(\gamma'(s), \gamma'(s))$, and $l=3\tau/4 + s \tau/4$.

This inequality together with the decay property of $R$ and
compatibility of distances to conclude
\[
G(x, \tau; x_0, \tau/2) \ge \frac{c(\kappa, D_0)}{ |t|^{n/2}} e^{-b_1 d(x, x_0, t)^2/\tau}.
\]This finishes the proof of the theorem. \qed

\section{Applications to ancient solutions and the Poincar\'e conjecture}

In this section we use Theorem \ref{th1.1} to give a different proof for
Perelman's classification result of backward limits of ancient $\kappa$
solutions.

\begin{theorem} (Perelman)
\lab{thbacklimit}
Let $g(\cdot, t)$ with $t \in (-\infty, 0]$ be a nonflat, $3$ dimensional ancient
$\kappa$ solution for some $\kappa>0$.  Then there exist sequences of points
$\{q_k \} \subset \M$ and times $t_k \to -\infty$, $k=1, 2, ...$, such that
the scaled metrics $g_k(x, s) \equiv R(q_k, t_k) g(x, t_k + s R^{-1}(q_k, t_k))$
around $q_k$
converge to a nonflat gradient shrinking soliton in $C^\infty_{loc}$
topology.
\proof
\end{theorem}

We divide the proof into several cases.

Case 1 is when the section curvature is zero somewhere and $\M$ is noncompact. Then Hamilton's strong
maximum principle for tensors show that $\M=\M_2 \times \R^1$ where $\M_2$ is
a $2$ dimensional, nonflat ancient $\kappa$ solution.  According to Hamilton, $\M_2$
is either $S^2$ or $RP^2$. So the theorem is already proven in this case.
This case can also be covered in Case 4 below together.

Case 2 is when the section curvature is zero somewhere and $\M$ is compact.

Then, again using maximum principle,
 Hamilton (see Theorem 6.64 in [CLN] e.g) showed that $\M$ is the metric quotient of
 $\R^3$ with the flat metric or that of $S^2 \times \R^1$. So the theorem is also
 proven in this case.

 Case 3 is when the sectional curvature is positive everywhere and $\M$ is a
 type II ancient solution. i.e. $\sup_{t<0} |t| \ R(\cdot, t) = \infty$.

In this case Hamilton \cite{H:2} showed by a scaling argument and
his matrix maximum principle that the backward limit is a steady
gradient soliton. See also Theorem 9.29 in \cite{CLN:1}, in which a proof is given 
for the non-compact case.  However the compact case can be proven in the same way
with the $\kappa$ non-collasping assumption. So one can
take a scaling
 limit to a shrinking gradient soliton. See Theorem 9.66 in \cite{CLN:1} e.g.
 If the ancient solution arises from the blow up of finite time type II
 singularity,
 then Hamilton \cite{H:2} even proved that $\M$ is a steady gradient soliton.
 If $\M$ is compact, then it is well known that $\M$ is an Einstein
 manifold.  Since the curvature is positive, $\M$ has to be $S^3$.

 So there is only one case left. 
 
 Case 4:   $\M$ has positive sectional curvature and
 is of type I ancient solution.

 If $\M$ is compact, N. Sesum already proved the theorem in this case \cite{S:1}.   Actually
 she proved a stronger result, namely, $\M$ is a shrinking gradient
 soliton. See also p 302 \cite{CZ:1} and
 the work of X.D. Cao \cite{Cx:1}.

 So we will assume that $\M$ is noncompact and of type I for the rest of the proof.
  In fact our proof works in both compact and noncompact cases.

  By the $k$ noncollapsing assumption and the bound $R(\cdot, t)
 \le \frac{D_0}{1+|t|}$, we can find a sequence $\tau_k \to \infty$ such that the following holds:

 the pointed manifolds $({\M}, g_k, y_k)$ with the metric
 \[
 g_k \equiv \tau^{-1}_k g(\cdot, - s \tau_k)
 \] converge, in $C^\infty_0$ sense, to a pointed manifold $({\M}_\infty, g_\infty(\cdot, s), y_\infty)$. Here $s>0$.
 

 We aim to prove that $g_\infty$ is a gradient, shrinking Ricci soliton.  Note that we are
 scaling by $\tau^{-1}_k$. By the upper and lower bound on the scalar curvature, this
 scaling is equivalent to scaling by the scalar curvatures.
 We define, for $x \in \M$ and $s \ge 1$, the functions
 \[
 u_k=u_k(x, s) \equiv \tau^{n/2}_k \ G(x, s \tau_k; x_0, 0).
 \]Here $G$ is the heat kernel of the conjugate heat equation and $x_0$ is a fixed point.  We choose $y_k=x_0$ in the scaled of metrics above.  By Theorem \ref{th1.1}
 (actually (\ref{Gdbound})), we know that $u_k(x, s) \le U_0$ uniformly for all $k=1, 2, ...$, $x \in \M$ and
 $s$ in a compact interval. Here $U_0$ is a positive constant.  Note that $u_k$ is a
 positive solution of the conjugate heat equation under the metric on $(\M, g_k(s))$ i.e.
 \[
 \Delta_{g_k} u_k - R_{g_k} u_k - \partial_s u_k  =0.
 \]We have seen that
 $u_k$ and $R_{g_k}$ are uniformly bounded on compact intervals of $s$ in $(0, \infty)$,
and also the Ricci curvature is nonnegative and the curvature
tensors are uniformly bounded. The standard parabolic theory shows
that $u_k$ is H\"older continuous uniformly with respect to $g_k$.
Hence we can extract a subsequence, still called $\{
 u_k \}$, which converges in $C^\alpha_{loc}$ sense, modulo diffeomorphism, to
 a $C^\alpha_{loc}$ function $u_\infty$ on $({\M}_\infty, g_\infty(s), y_\infty)$.

 Using integration by parts, it is easy to see that $u_\infty$ is a weak solution of
 the conjugate heat equation on  $({\M}_\infty, g_\infty(s))$, i.e.
 \[
 \int \int \left( u_\infty \Delta \phi - R u_\infty \phi + u_\infty \partial_s \phi \right)
 d g_\infty(s) ds =0
 \]for all $\phi \in C^\infty_0({\M}_\infty \times (-\infty, 0])$.

 By standard parabolic theory, the function $u_\infty$, being bounded on compact
 time intervals,
 is a smooth solution of the conjugate heat equation on $({\M}_\infty, g_\infty(s), y_\infty)$.  We need to show that $u_\infty$ is not zero.  One can even show that it is actually the fundamental solution of the conjugate 
 heat equation with pole at $y_{\infty}$   (the image of the same $x_0$ in the limiting 
 manifold).    Let $u=u(x, \tau)=G(x, \tau; x_0, 0)$. We claim that for 
 a constant $a>0$ and all $\tau \ge 1$,
 \[
 u(x_0, \tau) \ge \frac{a}{\tau}.
 \]Here is the proof.  Define $f$ by 
 \[
 (4 \pi \tau)^{-n/2} e^{-f} = u.
 \]By Corollary 9.4 in \cite{P:1}, which is a consequence of his differential Harnack inequality for fundamental solutions,  we have, for $\tau=-t$,
 \[
 -\partial_t f(x_0, t) \le \frac{1}{2} R(x_0, t) - \frac{1}{2 \tau} f(x_0, t).
 \]Since $R(x_0, t) \le c/\tau$, we can integrate the above from $\tau=1$ to get
 \[
 f(x_0, \tau) \le c + \frac{f(x_0, 1)}{\tau} \le C.
 \]Here we have use the fact that $f(x_0, 1)$ is bounded, by the standard short time bounds for $G=G(x_0, 1; x_0, 0)$.  This proves the claim.  By definition of $u_k$ as a scaling of $u$,  we know that $u_k(x_0, s) \ge b>0$ for $s \in [1, 4]$.  Here $b$ is independent of $k$. Therefore $u_\infty(x_0, s) \ge b>0$.
The maximum principle shows $u_\infty$ is positive everywhere.

 Let us recall that Perelman's $W$ entropy for each $u_k$ is
 \[
 W_k(s) = W(g_k, u_k, s) =
 \int \left[ s ( |\nabla f_k|^2 + R_k) +f_k - n \right] u_k dg_k(s)
 \]where $f_k$ is determined by the relation
 \[
 (4 \pi s)^{-n/2} e^{-f_k} = u_k;
 \]and $R_k$ is the scalar curvature under $g_k$.
 By the uniform upper bound for $u_k$, we know that there exist $c_0>0$ such that
 \[
 f_k = - \ln u_k - \frac{n}{2} \ln ( 4 \pi s) \ge - c_0
 \]for all $k=1, 2, ...$ and $s \in [1, 3]$.  Here the choice of this interval for $s$
 is just for convenience. Any finite time interval also works.
 Since $\M$ is noncompact, one needs to justify the integral in $W_k(s)$ is finite.
 For fixed $k$, $u_k$ has a generic Gaussian upper and lower bound with coefficients depending on $\tau_k$ and curvature tensor and their derivatives, as shown in \cite{G:1}.  The manifold has nonnegative Ricci 
 curvature and bounded curvature.   So the term $f_k u_k$ which is essentially $-u_k \ln u_k$ is integrable. The term $|\nabla f_k|^2 u_k = |\nabla u_k|^2 /u_k$ which is integrable by Theorem 3.3 in \cite{Z:1}, given in the appendix.  These together imply that $W_k(s)$ is well defined.

 Since $\int_{\M} u_k dg_k =1$, we know that
 \be
 \lab{wklowerb}
 W_k(s) \ge -c_0-n
 \ee  for all $k=1, 2...$ and $s \in [1, 3]$.

 There is an alternative proof of the lower bound for $W_k$.
 Actually $W_k(s)$ is uniformly bounded from below if $u_k$ is replaced by any $v\in
 W^{1, 2}$ such that $\Vert v \Vert_2 =1$.  This can be seen since $({\M}, g_k(s), y_k)$,
 $s \in [1, 3]$ has uniformly bounded curvature operator and are $\kappa$ noncollapsed.
 Therefore, a uniform Sobolev inequality holds, which implies the lower bound
 of $W_k(s)$. The later is nothing but a lower bound on the best constants of log
 Sobolev inequalities.

 By scaling it is easy to see that
 \[
 W_k(s) = W(g, u, s \tau_k),
 \]where $u=u(x, l) =G(x, l, x_0, 0)$. According to \cite{P:1},
 \be
 \lab{dwds}
 \frac{d W_k(s)}{ds} =
 -  2 s \int | Ric_{g_k} + Hess_{g_k} f_k - \frac{1}{2s}
 g_k |^2 u_k dg_k(s) \le 0.
 \ee  Note that the integral on the right hand side is finite by a similar argument as in 
 the case of $W_k(s)$.  So, for a fixed $s$, $W_k(s)=W(g, u, s \tau_k)$ is a non-increasing function of $k$.
 Using the lower bound on $W_k(s)$  (\ref{wklowerb}), we can find a function
 $W_\infty(s)$ such that
 \[
 \lim_{k \to \infty} W_k(s) = \lim_{k \to \infty} W(g, u, s \tau_k) = W_\infty(s).
 \]

 Now we pick $s_0 \in [1, 2]$. Clearly we can find a subsequence $\{ \tau_{n_k} \}$, tending to infinity,
 such that
 \[
 W(g, u, s_0 \tau_{n_k})
 \ge W(g, u, (s_0+1) \tau_{n_k}) \ge W(g, u, s_0 \tau_{n_{k+1}}).
 \]Since
\[
\lim_{k \to \infty} W(g, u, s_0 \tau_{n_k})  =
\lim_{k \to \infty} W(g, u, s_0 \tau_{n_{k+1}})  = W_\infty(s_0),
\]we know that
\[
\lim_{k \to \infty} [ W(g, u, s_0 \tau_{n_k}) - W(g, u, (s_0+1) \tau_{n_k}) ] =0.
\]That is
\[
\lim_{k \to \infty} [ W_{n_k}(s_0) - W_{n_k}(s_0+1) ] =0.
\]Integrating (\ref{dwds}) from $s_0$ to $s_0 + 1$, we use the above to conclude that
\[
 \lim_{k \to \infty} \int^{s_0+1}_{s_0} \int s | Ric_{g_{n_k}} + Hess_{g_{n_k}} f_{n_k} - \frac{1}{2s}
 g_{n_k } |^2 u_{n_k} dg_{n_k}(s) ds = 0.
\]Therefore
\[
Ric_\infty + Hess_\infty f_\infty - \frac{1}{2s} g_\infty=0.
\]So the backward limit is a gradient shrinking Ricci soliton. 

Finally we need to show the soliton is non-flat.  We can assume the  original ancient solution is not a gradient shrinking soliton. Otherwise there is nothing to prove.  Hence, we know that $W_k(s) < W_k(0) =W_0=0$ where $W_0$ is the Euclidean $W$ 
entropy with respect to the standard Gaussian.  Hence $W_\infty(s) \le W_k(s) <W_0$.
If the gradient shrinking soliton $g_\infty$ were flat, it is known to be $\R^3$. Hence $W_\infty(s)=W_0$, a contradiction.
\qed
\medskip

{\it Remark.   Case 4 with positive curvature tensor can also be dealt with 
by the method in \cite{CL:1}. There Chow and Lu actually constructed an embedded region of the flow , which is close to $S^2 \times \R$. They even 
do not need to assume the soliton is $\kappa$ non-collapsed on all scales.
In fact, there does not exist type I, noncompact, $\kappa$ ancient solution with positive curvature tensor, after all.  This is due to Perelman's classification of backward limits.

Also the on diagonal lower bound of the fundamental solution $G$ in the middle of the proof can be extended to full lower bound by the theorem in the
appendix.  But we do not need it here. }

\medskip

In the last part of the section, we discuss the ramification of the
above method to the proof of the Poincar\'e conjecture.  After the
classification of the backward limits and $\kappa$ non-collapsing
with surgeries, the only part of Perelman's proof of the Poincar\'e
conjecture that requires the reduced distance and volume is the
universal non-collapsing of ancient $\kappa$ solutions.
Interestingly, a different proof of this fact already exists in
Section 3.2 of the paper of Chen and Zhu \cite{ChZ:1}, where certain
more general $4$ dimensional result is proven (see the paragraph
after the proof of Proposition 3.4 there). In the $3$ dimensional
case, the proof looks longer than Perelman's original proof. However
it is basically a reshuffling of certain arguments suggested by
Perelman, all which are needed to prove the canonical neighborhood
property for ancient $\kappa$ solutions.  In this sense, the proof
of the universal noncollapsing is a by product of canonical
neighborhood property for ancient $\kappa$ solutions.  Indeed, the
canonical neighborhood property for ancient $\kappa$ solutions can
be proven exactly the same way without the universal noncollapsing
property, except that the constants in the property depend on the
noncollapsing constant $\kappa$.  But this is enough to show that
after a conformal change of metric using the scalar curvature
function, the ancient solution is $\epsilon$ close to model
manifolds which are universal noncollapsed. Therefore the former is
also universal noncollapsed.

Let us state the result and sketch the proof.

\begin{proposition} (Perelman)
\lab{prk0noncollapsing}
There exists a positive constant $\kappa_0$ with the following property.
 Suppose we have a non-flat, $3$ dimensional ancient $\kappa$
 solution arising from finite time singularity of a Ricci flow,
  for some $\kappa>0$.
Then either the solution is $\kappa_0$ non-collapsed on all scales
or it is a metric quotient of the round $3$ sphere.
 \proof (
sketched as a special case of Chen and Zhu's proof in Section 3.2, the
statement after Proposition 3.4 \cite{ChZ:1})
\end{proposition}

Note we use an extra assumption that $\kappa$
 solution is arising from finite time singularity of a Ricci flow.
 This will make the proof more transparent since type II $\kappa$
 solution in this case is just steady gradient Ricci soliton as proven by
 Hamilton \cite{H:2}.

If the three dimensional $\M$ is compact, then they are explicitly
known to be gradient solitons as mentioned in Cases 1-4 in the proof
of the previous theorem. Anyway they are not needed in singularity
analysis leading to the Poincar\'e conjecture. So we just need to
prove that noncompact $3$ dimensional $\kappa$ ancient solutions are
universal non-collapsed on all scales.  The proof is divided into
$3$ steps.

{\bf step 1.}  one proves the compactness of ancient $\kappa$
solutions with any fixed $\kappa>0$. i.e.

{\it The set of nonflat $3$ dimensional ancient $\kappa$ solutions,
for any fixed $\kappa>0$,  is compact modulo scaling in the
following sense: for any sequence of such solutions and marking
points in space time  $(x_k, 0)$ with $R(x_k, 0)=1$, one can extract
a $C^\infty_{loc}$ converging subsequence whose limit is also an
ancient $\kappa$ solution.}

The proof is identical to that in \cite{P:1}, the Theorem in Section
11.7. Note that no universal non-collapsing is needed here. This
actually is the original order of proof by Perelman.

{\bf step 2.}  One proves certain elliptic type estimates for the scalar curvature.

{\it  There exist a positive constant $\eta$ and a positive
increasing function $w: [0, \infty) \to (0, \infty)$ with the
following property.  Let $(\M, g_{ij}(t))$, $-\infty < t \le 0$ is a
$3$ dimensional ancient $\kappa$ solution for a fixed $\kappa>0$.
Then

(i) for every $x, y \in \M$ and $t \in (-\infty, 0]$, there holds
\[
R(x, t) \le R(y, t) \, w(R(y, t) d^2(x, y, t));
\]

(ii) for all $x \in \M$ and $t \in (-\infty, 0]$, there hold
\[
|\nabla R| \le \eta R^{3/2}(x, t), \qquad |\partial_t R|(x, t) \le
\eta R^2(x, t).
\]

(iii)  Suppose for some $(y, t_0)$ in space time and a constant
$\zeta>0$ there holds
\[
\frac{| B(y, R(y, t_0)^{-1/2}, t_0)|_{t_0}}{R(y, t_0)^{-3/2}} \ge
\zeta.
\] Then there exist a positive functioins $w$ depending only on
$\zeta$ such that, for all $x \in \M$,
\[
R(x, t_0) \le R(y, t_0) \, w(R(y, t_0) d^2(x, y, t_0)).
\]
}

The proof of statements (i) and (ii) is almost a carbon copy of
Theorem 6.4.3 in \cite{CZ:1} (3 d case) or Proposition 3.3 (4 d
case) in \cite{ChZ:1}, or the corresponding results in \cite{KL:1}
and \cite{MT:1}. The one difference is that one uses $\kappa$
non-collapsing assumption instead the universal non-collapsing that
is being proved. Therefore the constant $\eta$ and the function $w$
may depend on $\kappa$. Part (iii) is the remark after Proposition
3.3 (4 d case) in \cite{ChZ:1}, which includes the $3$ dimension
case as a special situation. Its proof is a moderate refinement  of
that of  statement (i), by keeping a careful track of constants.

{\bf step 3.}  For any point $(x, t)$, one shows that either it is a
center of the $\epsilon$ neck, or it lies in a compact manifold with
boundary, called $\M_\epsilon$.  After scaling by scalar curvature
at one of its boundary points, this manifold is $\epsilon$ close to
a compact manifold of finite diameter and whose scalar curvature is
bounded between two positive constants which are {\it independent}
of the noncollapsing constant $\kappa$. This step follows
Proposition 3.4 in \cite{ChZ:1} which is a 4 dimensional result that
includes the 3 dimension one as a special case. They use a blow up
argument, taking advantage of the property that a boundary point of
$\M_\epsilon$ is the centered of a ball which is $2\epsilon$ close
to that of $S^2 \times \R$ after scaling.  Then they use (iii) in
step 2 to obtain the bounds on scalar curvature.  The bounds depend
only on the noncollapsing constant of $S^2 \times \R$.

This means that after scaling by scalar curvature, every point on
the ancient solution has a ball of fixed diameter that is $\epsilon$
close to a model manifold which is universal non-collapsed.
Therefore ancient $\kappa$ solution is also universal non-collapsed.
\qed
\medskip

Let us close by presenting the flow chart of a simplified
proof of the Poincar\'e
conjecture without reduced distance or volume.

Step 1.  $W$ entropy and its monotonicity  (\cite{P:1}).  See also
\cite{Cetc:1}, \cite{CZ:1}, \cite{KL:1}, \cite{MT:1}.

Step 2.  Local non-collapsing result via Step 1 (\cite{P:1}).
                See also
\cite{Cetc:1}, \cite{CZ:1}, \cite{KL:1}, \cite{MT:1}.

Step 3. getting ancient $\kappa$ solutions by blowing up of
singularity using Step 2 and Hamilton's  compactness theorem
(\cite{P:1}).  See also
 \cite{Cetc:1}, \cite{CZ:1}, \cite{KL:1}, \cite{MT:1}.

Step 4.  (i)  showing the backward limits of ancient $\kappa$
solutions are gradient shrinking solitons.  Earlier work of Hamilton
\cite{H:2} for type II case and \cite{CL:1} or this paper for type I case.

  (ii)  universal non-collapsing of ancient $\kappa$ solutions.  Section 3.2 of
  \cite{ChZ:1}.

  (iii) curvature and volume estimates for ancient solutions (\cite{P:1}).  See also
\cite{Cetc:1}, \cite{CZ:1}, \cite{KL:1}, \cite{MT:1}.

Step 5.  classification of  gradient shrinking solitons.
\cite{P:1}. See also
 \cite{Cetc:1}, \cite{CZ:1}, \cite{KL:1}, \cite{MT:1}.

Step 6.  canonical neighborhood property \cite{P:1}.  That is:
regions of high scalar curvature resemble the ancient solution after
appropriate scaling.

See also
\cite{Cetc:1}, \cite{CZ:1}, \cite{KL:1}, \cite{MT:1}.

Step 7.  surgery procedure, including properties of the standard solution
\cite{P:2}.  See also
\cite{CZ:1}, \cite{KL:1}, \cite{MT:1}.

Step 8.  local $\kappa$ non-collapsing with surgeries \cite{Z:2}.

Step 9.  canonical neighborhood property with surgeries \cite{P:2}.
 See also
\cite{CZ:1}, \cite{KL:1}, \cite{MT:1}.

Step 10. existence of Ricci flow with surgeries, i.e. proving there
are finitely many surgeries within finite time. \cite{P:2}.
 See also
\cite{CZ:1}, \cite{KL:1}, \cite{MT:1}.

Step 11.  Finite time extinction of Ricci flow on simply connected manifolds
\cite{P:3}.  See also \cite{CM:1} and \cite{MT:1}.

\section{Appendix}

Here we state and prove Theorem 3.3 in \cite{Z:1}, which was used at
the end of Step 2 in the proof of Theorem \ref{th1.1}. See also
\cite{CH:1}.

\begin{theorem}
Let ${\bf M}$ be a compact or complete noncompact Riemannian
manifold with bounded curvature and equipped with a family of
Riemannian metric evolving under the forward Ricci flow $\partial_t
g = - 2 Ric$ with $t \in [0, T]$. Suppose $u$ is any positive
solution to $\Delta u - \partial_t u =0$ in ${\bf M} \times [0, T]$.
Then, it holds
\[
\frac{|\nabla u(x, t) |}{u(x, t)} \le  \sqrt{ \frac{1}{t}}
    \sqrt{\log \frac{M}{u(x, t)}}
\]for $M = \sup_{{\bf M} \times [0, T]} u$ and $(x, t) \in {\bf M} \times [0,
T]$.

Moreover,  the following interpolation inequality holds for any
$\delta>0$, $x, y \in {\bf M}$ and $0< t \le T$:
\[
u(y, t) \le c_1 u(x, t)^{1/(1+\delta)} M^{\delta/(1+\delta)} e^{c_2
d(x, y, t)^2/t}.
\]Here $c_1, c_2$ are positive constants depending only on
$\delta$.
\end{theorem}
\medskip

{\bf Proof}

This is almost the same as that of Theorem 1.1 in \cite{H:1}. By
direct calculation
\[
\Delta (u \log \frac{M}{u}) - \partial_t (u \log \frac{M}{u})
=-\frac{|\nabla u|^2}{u},
\]
\[
(\Delta -\partial_t)(\frac{|\nabla u|^2}{u}) = \frac{2}{u}
\bigg{|}\partial_i \partial_j u - \frac{\partial_i u \partial_j
u}{u} \bigg{|}^2 \\
 \ge 0.
\]The first inequality follows immediately from the maximum principle since
\[
t \frac{|\nabla u|^2}{u} - u \log \frac{M}{u}
\]is a sub-solution of the heat equation.

To prove the second inequality, we set
\[
l(x, t) = \log (M/u(x, t)).
\]Then the first inequality implies
\[
|\nabla \sqrt{l(x, t)} | \le 1/\sqrt{t}.
\]Fixing two points $x$ and $y$, we can integrate along a geodesic to
reach
\[
\sqrt{\log (M/u(x, t))} \le \sqrt{\log (M/u(y, t))} + \frac{d(x, y,
t)}{\sqrt{t}}.
\]The result follows by squaring both sides.
\qed


\medskip

\noindent e-mail:  qizhang@math.ucr.edu

date:  December 2008

\enddocument